\documentclass[a4paper]{article}

\usepackage{amsmath}
\usepackage{amssymb}
\usepackage{amsthm}
\usepackage[matrix,arrow]{xy}
\usepackage{xspace}

\newcommand{\CP}{\ensuremath {\mathbb{CP}}\xspace}
\newcommand{\N}{\ensuremath {\mathbb{N}}\xspace}
\newcommand{\F}{\ensuremath {\mathbb{F}}\xspace}
\newcommand{\Z}{\ensuremath {\mathbb{Z}}\xspace}
\DeclareMathOperator{\coker}{coker}
\DeclareMathOperator{\im}{im}
\newtheorem{theorem}{Theorem}[section]
\newtheorem{lemma}[theorem]{Lemma}
\newtheorem{proposition}[theorem]{Proposition}
\newtheorem{defn}[theorem]{Definition}
\newtheorem{corollary}[theorem]{Corollary}
\newtheorem{thm}{Theorem}
\newtheorem{cor}[thm]{Corollary}
\newtheorem{examples}[theorem]{Examples}

\newcommand{\m}[1]{\ensuremath {\mathcal{#1}}\xspace}
\newcommand{\mf}[1]{\ensuremath {\mathfrak{#1}}\xspace}
\newcommand{\smsh}{\land}
\newcommand{\abs}[1]{\left\lvert #1 \right\rvert}

\newcommand{\ldb}{[\![}
\newcommand{\rdb}{]\!]}
\newcommand{\lb}{[}
\newcommand{\rb}{]}
\newcommand{\noproof}{{\hspace*{\fill}\qed}}

\newcommand{\kn}{{K(n)}}
\newcommand{\plusf}[1][F]{+_{\raisebox{-2pt}{\(\scriptscriptstyle \!\!#1\)}}}
\newcommand{\chcl}[1][F]{x^{\scriptscriptstyle #1}}
\newcommand{\pseries}[2]{\lb #1 \rb_{\scriptscriptstyle #2}}
\newcommand{\bp}{B P}

\newcommand{\rc}[1]{\widetilde{#1}}
\newcommand{\spc}[2][k]{\underline{#2}_{#1}}

\newcommand{\rkn}{\widetilde{K}(n)}

\newcommand{\hyp}{-\hspace{0pt}}

\DeclareMathOperator{\colim}{colim}
\DeclareMathOperator{\map}{Map}
\DeclareMathOperator*{\bstar}{\bigstar}

\newcommand{\co}{\colon\thinspace}

\numberwithin{equation}{section}

\begin{document}

\title{%
Stable and Unstable Operations in mod \(p\)\\Cohomology Theories%
}

\author{Andrew Stacey\\%
Sarah Whitehouse}

\maketitle

\begin{abstract}
We consider operations between two multiplicative, complex orientable
cohomology theories. Under suitable hypotheses, we construct a map
from unstable to stable operations, left-inverse to the usual map from
stable to unstable operations.  In the main example, where the target
theory is one of the Morava K-theories, this provides a simple and
explicit description of a splitting arising from the Bousfield-Kuhn
functor.
\end{abstract}

\section{Introduction}
\label{sec:intro}

Given two graded cohomology theories, \({E}^*(-)\) and
\({F}^*(-)\), we can consider various types of operations from one
to the other.  There are two main types: stable and unstable; and
within the unstable operations are the additive operations.

These are simplest to describe in categorical language.  A cohomology
theory is a functor satisfying certain properties.  At various levels
of forgetfulness we have the following functors:
\begin{enumerate}
\item A functor \({E}^*(-)\) from the (homotopy) category of based
topological spaces to the category of graded abelian groups which
intertwines the two suspension functors.

\item A sequence of functors \(({E}^k(-))_{k \in \Z}\) from the
(homotopy) category of based topological spaces to the category of
abelian groups.

\item A sequence of functors \(({E}_U^k(-))_{k \in \Z}\) from the
(homotopy) category of based topological spaces to the category of
sets.
\end{enumerate}

The three types of operation from \({F}^*(-)\), the source theory, to
\({E}^*(-)\), the target theory, are:
\begin{enumerate}
\item Stable operations: for \(l \in \Z\), \(\m{S}^l(F,E)\) is the set
of natural transformations \(r \co {F}^*(-) \to {E}^*(-)\) of degree \(l\).

\item Additive operations: for \(k,l \in \Z\), \(\m{A}_k^{k+l}(F,E)\)
is the set of natural transformations \(r_k \co {F}^k(-) \to {E}^{k+l}(-)\).

\item Unstable operations: for \(k, l \in \Z\), \(\m{U}_k^{k+l}(F,E)\)
is the set of natural transformations \(r_k \co {F}_U^k(-) \to
{E}_U^{k+l}(-)\).
\end{enumerate}

There is an obvious restriction map \(\m{S}^l(F,E) \to
\m{U}_k^{k+l}(F,E)\) for each \(k,l \in \Z\).  In brief, our main
theorem shows this map has a left-inverse under certain conditions on
\({E}^*(-)\) and \({F}^*(-)\).  The full statement of the
theorem is as follows.

\begin{thm}
\label{th:main}
Let \({E}^*(-)\) and \({F}^*(-)\) be two multiplicative graded
cohomology theories which are commutative and complex orientable.  Let
\(E^* := {E}^*(\text{pt})\) be the coefficient ring of \({E}^*(-)\).  We
assume that the following conditions hold.
\begin{enumerate}
\item \(E^*\) has characteristic \(p\).

\item The formal group law of \({E}^*(-)\) has finite height, say \(n\).

\item The coefficient of the first term in the \(p\)\hyp{}series for
\({E}^*(-)\) is invertible.

\item The various \(E^*\)\hyp{}modules of operations from
\({F}^*(-)\) to \({E}^*(-)\) are the duals over \(E^*\) to the
corresponding \(E^*\)\hyp{}modules of co-operations.
\end{enumerate}

Under these conditions, for each \(k,l \in \Z\) there is a map
\[
\Delta^\infty \co \m{U}_k^{k+l}(F,E) \to \m{S}^l(F,E)
\]
which is left-inverse to the restriction map.
\end{thm}

We postpone to the next section an explanation of what all the
conditions mean.  The map \(\Delta^\infty\) has several pleasant
properties; to describe most of these we need to know more about the
structure of the spaces of the various types of operation, knowledge
that we also postpone for the next section.

The map itself has a very simple description.  To give this in its
most topological form we recall that operations between cohomology
theories are closely related to homotopy classes of maps between
certain spaces and between certain spectra associated to the
cohomology theories.  In this language, the restriction map from
stable operations to unstable operations is nothing more than the
infinite-loop space functor, \(\Omega^\infty\).  Thus we obtain the
following corollary of theorem~\ref{th:main}.

\begin{cor}
\label{cor:maps}
Let \({E}^*(-)\) and \({F}^*(-)\) be cohomology theories as in
theorem~\ref{th:main}.  Let \(E\) and \(F\) be representing spectra.
Let \(\Omega^\infty\) denote the infinite-loop space functor from
spectra to topological spaces.  Then there is a map:
\[
  \Delta^\infty \co \lb \Omega^\infty F, \Omega^\infty E \rb_+ \to \{F,
    E\}^0
\]
left-inverse to the map induced by the \(\Omega^\infty\)\hyp{}functor.
\end{cor}

The subscript adorning \(\lb X, Y\rb_+\) is to denote homotopy classes
of maps which preserve the basepoint.

Composing \(\Delta^\infty\) with the map coming from the
\(\Omega^\infty\)\hyp{}functor we produce a projection on \(\lb
\Omega^\infty F, \Omega^\infty E\rb_+\) with the property that a
homotopy class lies in the image of this projection if and only if it
is an infinite loop map and, moreover, the delooping of this map is
unique.

This projection is easy to describe.  There are certain maps:
\begin{align*}
  v^E_n &\co \Omega^\infty E \to \Omega^{2(p^n-1)} \Omega^\infty E \\
  v^F_n &\co \Omega^\infty F \to \Omega^{2(p^n-1)} \Omega^\infty F
\end{align*}
which come from the \(p\)\hyp{}series of the formal group law associated to
each cohomology theory.  The conditions on \({E}^*(-)\) guarantee that
\(v_n\) is invertible.  The projection is:
\[
  \rho \mapsto (v^E_n)^{-1} \left( \Omega^{2(p^n-1)} \rho \right) v^F_n.
\]

The conditions in the theorem are really all about the target theory,
\({E}^*(-)\), even the last one.  The main examples to which we wish
to apply this theorem are where \({E}^*(-)\) is one of the Morava
K\hyp{}theories, \(\kn^*(-)\), at a prime \(p\).  We discard the case
\(n = 0\) as that is just rational cohomology and we take an odd prime
to ensure that the multiplication is commutative.  With this choice
for the target theory there is no restriction on choosing \({F}^*(-)\)
as the four conditions in theorem~\ref{th:main} are automatically
satisfied.  In this case, corollary~\ref{cor:maps} is an elaboration
of an application of the Bousfield-Kuhn functor.

This functor, written \(\Phi_n\), goes from the homotopy category of
based \(p\)\hyp{}local spaces to the homotopy category of
\(p\)\hyp{}local spectra.  Its key property is that if \(G\) is a
\(p\)\hyp{}local spectrum then \(\Phi_n \Omega^\infty (G)\) is
\(L_{\kn} G\), the \(\kn\)\hyp{}localisation of \(G\).  In particular,
if \(G\) is already \(\kn\)\hyp{}local then \(\Phi_n \Omega^\infty(G)
\simeq G\).  This is the case for \(\kn\) itself.  Thus the functorial
properties of \(\Phi_n\) yield a map:
\[
   [\Omega^\infty G, \Omega^\infty \kn]_+ \to \{L_{\kn} G, \kn\}^0.
\]
One of the defining properties of the \(\kn\)\hyp{}localisation is that
there is a natural isomorphism \(\{L_{\kn} G, \kn\} \cong \{G,
  \kn\}\).  Therefore we have a map
\[
  \Theta_n \co [\Omega^\infty G, \Omega^\infty \kn]_+ \to \{G, \kn\}^0.
\]
By a similar device we can remove the requirement that \(G\) be
\(p\)\hyp{}local.  Then \(\Theta_n\) can be compared to the map
\(\Delta^\infty\) in corollary~\ref{cor:maps}.

\begin{thm}
\label{th:bk}
Let the source theory be a complex orientable, graded, commutative
multiplicative cohomology theory.
Then, with target theory \(\kn^*(-)\), \(\Delta^\infty = \Theta_n\).
\end{thm}

This paper is structured as follows.  In the next section we describe
the features of cohomology theories that we need.  In
section~\ref{sec:pseries} we look at the \(p\)\hyp{}series coming from the
complex orientation and use this to define certain key co-operations.
These are the essential ingredients of the proof of
theorem~\ref{th:main}.  In section~\ref{sec:splitting} we prove our
main technical result, proposition~\ref{prop:desusp}, which involves
the relationships between the spaces of co-operations.  It is then a
short step to our main result, theorem~\ref{th:opsmaps} in
section~\ref{sec:maps}, which is a more detailed version of 
theorem~\ref{th:main} and of corollary~\ref{cor:maps}.  We conclude by
proving theorem~\ref{th:bk} in section~\ref{sec:bk}.

There is considerable detail in the papers \cite{jb4} and
\cite{jbdjww} about operations in cohomology theories.  Most of the
background that we need can be found in those papers.  The papers
\cite{drww} and \cite{ww} are the original sources for some of the
structure that we use.

There is some overlap in our main theorem with the work of
\cite{tknspt}.  The first splitting of \(K_*(\bp)\) in that paper is
dual to our splitting.  We do not go on to consider further
splittings, as is done in \cite{tknspt}, because the first
splitting has a good topological description which is missing in the
higher ones.  The proof of our theorem and that of \cite{tknspt} run
along similar lines.

Related work using the Bousfield-Kuhn functor has appeared in
\cite{ab6}, \cite{nk3}, and \cite{math.AT/0407022}.

Finally, we note some conventions that we shall use throughout this
paper.  Firstly, we work throughout in the \emph{homotopy} categories
of spaces and spectra and shall use the short-hand ``map'' for a
morphism in the appropriate category.  Thus what we mean when we say
``map'' is really a homotopy class of maps in the conventional sense.
We trust that this will not be overly confusing.

Secondly, following \cite{jb4} and \cite{jbdjww} we grade homology
negatively.  In order to get the pairing between homology and
cohomology correct one theory has to be graded negatively.  As in
\cite{jb4} and \cite{jbdjww}, for us the cohomology theory is the
object of study whereas the homology theory is a tool we shall use in
the analysis.

Thirdly, and unlike \cite{jb4} and \cite{jbdjww}, we shall always be
careful to distinguish between \emph{spaces} and \emph{spectra}.  The
convention of \cite{jb4} and \cite{jbdjww} is to use the same notation
for a space and its suspension spectrum.  This is a convenient
shorthand but as our paper is all about the passage from spaces to
spectra it is a shorthand we feel morally obliged to do without.

Fourthly, we shall need to work with both based and unbased spaces.
We shall distinguish between morphisms in the two categories with the
notations \(\lb X, Y \rb\) for homotopy classes of \emph{all} maps and
\(\lb X, Y \rb_+\) for homotopy classes of \emph{based} maps.  We
recall that when the target, \(Y\), is an \(H\)\hyp{}space and the
source, \(X\), is a based space then there is a natural projection
\(\lb X, Y \rb \to \lb X, Y \rb_+\).

\section{Ingredients}
\label{sec:ingredients}

In this section we shall describe the various ingredients needed for
our work.  This is not intended to be a detailed reference on
cohomology theories, rather our aim is to establish our notation
whilst giving just enough information to allow the casual reader to
follow our argument without constantly referring to other works.  The
bulk of this can be found in the expository parts of \cite{jb4} and
\cite{jbdjww} and we largely follow their conventions.  The reader
familiar with \cite{jb4} and \cite{jbdjww} may wish to skip to the
next section.

\subsection{Generalised Cohomology Theories}
\label{sec:gencoh}

Let \({E}^*(-)\) be a multiplicative graded generalised cohomology
theory that is commutative and complex orientable.  Much of what we
are about to say applies to more general theories but as we shall only
use such theories we specialise at the outset.

As this is a multiplicative theory, the cohomology of a point is a
graded ring called the \emph{coefficient} ring.  We write this as
\(E^*\).

\paragraph{Representing Spaces and Spectrum.}

Brown's representability theorem, and its consequences, provide us
with a sequence of \(H\)\hyp{}spaces, \((\spc{E})_{k \in \Z}\), which
represent this cohomology theory.  That is, we have universal elements
\(\iota_k \in {E}^k(\spc{E})\) such that for any space \(X\) the map
\(\alpha \mapsto \alpha^* \iota_k\) is an isomorphism of abelian groups:
\[
\lb X, \spc{E} \rb \to {E}^k(X).
\]
The abelian group structure on the left-hand side comes from the
\(H\)\hyp{}space structure of \(\spc{E}\).  The universal class
\(\iota_k\) actually lies in the subgroup \(\rc{E}^k(\spc{E})\) and so
for any \emph{based} space \(X\) the above isomorphism identifies
\(\lb X, \spc{E} \rb_+\) with \(\rc{E}^k(X)\).

These spaces are unique up to equivalence.  It can be shown that the
suspension isomorphism of reduced cohomology, \(\rc{E}^k(X) \cong
\rc{E}^{k+1}(\Sigma X)\), defines an equivalence \(\spc{E} \to \Omega
\spc[k+1]{E}\).  These equivalences allow us to construct an
\(\Omega\)\hyp{}spectrum \(E\) from the \(\spc{E}\).  Using this
spectrum we can extend the cohomology theory to spectra by defining
\(\rc{E}^k(F) := \{F, E\}^k\) and define the associated homology
theory for both based spaces and spectra as \(\rc{E}_k(X) := \{S, E
  \smsh X\}^{-k}\).  This extends to unbased spaces by the usual
method of adding a disjoint basepoint: \(E_k(X) := \rc{E}_k(X_+)\).
Note that we are following the convention of \cite{jb4} in (redundantly)
writing the homology and cohomology of spectra as reduced.

In light of the fact that \(\rc{E}^k(F)\) and \(\{F, E\}^k\) are one
and the same for spectra, we make the same identification for spaces.
That is, we consider the isomorphism \(\lb X, \spc{E} \rb \cong
{E}^k(X)\) to be so natural as to be worth writing as an equality.
We shall still employ the language of both sides and talk of maps or
classes as best fits, but shall regard the two dialects as synonymous.

\paragraph{Structure Maps.}

All of the structure of the cohomology theory \({E}^*(-)\) is reflected
in the spectrum \(E\) and the spaces \(\spc{E}\).  Essentially, any
natural transformation of cohomology theories is represented by maps
between the associated spaces or spectra.  The existence of the map
can usually be deduced by applying the natural transformation to the
appropriate universal class.

As an example, we have the already-mentioned equivalence \(\spc{E}
\simeq \Omega \spc[k+1]{E}\) coming from the natural isomorphism
\(\rc{E}^k(X) \cong \rc{E}^{k+1}(\Sigma X)\).  To define the
associated map we apply the suspension isomorphism to the space
\(\spc{E}\):
\[
  \rc{E}^k(\spc{E}) \cong \rc{E}^{k+1}(\Sigma \spc{E}).
\]
By the representability theorem, the image of the universal class
\(\iota_k\) is represented by a (based) map \(\vartheta_k \co \Sigma
\spc{E} \to \spc[k+1]{E}\).  The naturality of the suspension
isomorphism implies that for a general space it is the composition:
\[
  \rc{E}^k(X) = \lb X, \spc{E} \rb_+ \xrightarrow{\Sigma} \lb \Sigma X,
  \Sigma \spc{E} \rb_+ \xrightarrow{{\vartheta_k}_*} \lb \Sigma X,
  \spc[k+1]{E} \rb_+ = \rc{E}^{k+1}(\Sigma X).
\]

In this fashion we deduce the existence of several maps which we now
list.

\begin{description}
\item[Suspension.] There is a map \(\vartheta_k \co \Sigma \spc{E} \to
\spc[k+1]{E}\) representing the suspension isomorphism \(\rc{E}^k(X) \cong
\rc{E}^{k+1}(\Sigma X)\).

\item[Stabilisation.] There is a map of spectra \(\sigma_k \co
\Sigma^\infty \spc{E} \to E\) of degree \(k\) representing the
isomorphism \(\rc{E}^k(X) \cong \rc{E}^k(\Sigma^\infty X)\).

\item[Multiplication.] There is a map of spectra, \(\phi \co E \smsh E
\to E\), of degree \(0\) and maps of spaces \(\phi_{k,l} \co \spc{E}
\smsh \spc[l]{E} \to \spc[k+l]{E}\) representing the multiplication in
the cohomology rings.

\item[Unit.] There is a map of spectra, \(\eta \co S \to E\), of degree
\(0\) and maps of spaces \(\eta_k \co S^k \to \spc{E}\) representing the
unit in the cohomology rings.
\end{description}

These maps satisfy various compatibility relations.  In particular,
the stable and unstable realms correspond under the stabilisation
maps.  We record one particular relation that will be of use later:
\begin{equation}
\label{eq:unit}
\vartheta_k = \phi_{1,k} (\eta_1 \smsh 1).
\end{equation}

Using the multiplication we can define the \emph{augmentation} maps.
The stable augmentation map is:
\[
  \epsilon_S \co \rc{E}_k(E) = \{S, E \smsh E\}^{-k}
  \xrightarrow{\phi_*} \{S, E\}^{-k} = E^{-k}.
\]
The unstable augmentations are:
\[
  \epsilon_k \co E_l(\spc{E}) \to \rc{E}_l(\spc{E}) \cong
  \rc{E}_l(\Sigma^\infty \spc{E}) \xrightarrow{{\sigma_k}_*}
  \rc{E}_{l-k}(E) \xrightarrow{\epsilon_S} E^{k-l}.
\]
We shall write \(\epsilon\) for \(\epsilon_k\) where we do not wish to
or cannot specify the index.

\paragraph{Duality.}

The augmentations define a pairing between cohomology and homology.
An element \(\alpha \in {E}^k(X)\) defines a push-forward in
homology for \(X\) a space or spectrum; respectively:
\begin{align*}
\alpha_* &\co {E}_l(X) \to {E}_l(\spc{E}), \\
\alpha_* &\co {E}_l(X) \to {E}_{l-k}(E)
\end{align*}
which we compose with the appropriate augmentation to end up in
\({E}^{k-l}\).

Under favourable circumstances the induced map \({E}^*(X) \to D
{E}_*(X)\) (the \(E^*\)\hyp{}dual of \({E}_*(X)\)) is an
isomorphism.  To truly understand this statement would require a
lengthy and, for our purposes, unnecessary discussion of the
topologies involved.  The precise circumstances are recorded in
\cite[theorem 4.14]{jb4}: if \({E}_*(X)\) is free as an
\(E^*\)\hyp{}module then \({E}^*(X)\) is the \(E^*\)\hyp{}dual of
\({E}_*(X)\).  When this occurs we shall say that \(X\) has
\emph{strong \(E\)\hyp{}duality}.  If this holds for all spaces and
spectra then we shall say that \({E}^*(-)\) has \emph{strong duality}.

\subsection{Operations and Co-operations}
\label{sec:opscoops}

Another piece of the baggage that comes with a generalised cohomology
theory is the family of operations.  As with all the other parts of
the structure of the cohomology theory, these are reflected in maps
between the representing spaces.  Also we can consider operations from
one cohomology theory to another.  Thus let \({F}^*(-)\) be another
generalised cohomology theory (also multiplicative, commutative, and
complex orientable).  We shall consider the operations from \({F}^*(-)\)
to \({E}^*(-)\).

\paragraph{Stable and Unstable Operations.}

As mentioned in the introduction, operations are simplest to describe
in the language of category theory.  In this setting, a cohomology
theory is a functor on the homotopy category of topological spaces and
an operation is simply a natural transformation between functors.
With a graded cohomology theory one has two types of operation
depending on whether one considers the cohomology theory as a whole,
leading to \emph{stable} operations, or one takes a single component
of it, leading to \emph{unstable} operations.  We allow degree shifts
in both cases.

In the stable case this description needs a little elaboration.
Considered as a whole, a cohomology theory is a functor between two
categories each of which has a suspension functor and the cohomology
theory intertwines these functors.  To qualify as a stable operation,
a natural transformation has also to respect the suspension functors.
Otherwise, using the restrictions mentioned below, a stable operation
would be simply a sequence of unstable operations with no relations
between the components.  Respecting the suspension functors imposes
some relations between successive components, as we will see in a
moment.

We label the set of stable operations of degree \(l\) from
\({F}^*(-)\) to \({E}^*(-)\) by \(\m{S}^l(F,E)\) and the set of
unstable operations from \({F}^k(-)\) to \({E}^l(-)\) by
\(\m{U}^l_k(F,E)\).  At the most basic level these are abelian groups
since operations take values in abelian groups.

There is an obvious way to define an unstable operation by restricting
a stable operation to a single component.  In this way, a stable
operation defines a sequence of unstable operations.  This suggests
the question as to whether a sequence of unstable operations patches
together to give a stable operation.  This will happen if the unstable
operations commute with the suspension isomorphisms, modulo a sign.
That is, suppose that for each \(k \in \Z\) we have an unstable
operation \(r_k \co {F}^k(-) \to {E}^{k+l}(-)\) then there is a stable
operation \(r \co {F}^*(-) \to {E}^{*+l}(-)\) restricting to \(r_k\)
(modulo the sign issue) if and only if each \(r_k\) maps reduced
cohomology to reduced cohomology and the following diagram commutes
for each space \(X\) up to the indicated sign:
\[
  \xymatrix{
    \rc{F}^k(X) \ar[r]^{\cong} \ar[d]^{r_k} \ar@{}[rd]|{(-1)^k}&
    \rc{F}^{k+1}(\Sigma X) \ar[d]^{r_{k+1}} \\
    \rc{E}^{k+l}(X) \ar[r]^{\cong} &
    \rc{E}^{k+l+1}(\Sigma X)
  }
  \]
The resulting stable operation need not be unique, however.  For that
one needs to know that a certain \(\lim^1\) term vanishes.  There are
technical conditions that guarantee this which, as we note later, hold
in our context.

A stable operation extends in the obvious way to an operation on the
cohomology of spectra.  There is no analogue of an unstable operation
in this case.

Using the same techniques as for the structure maps we can identify
operations with maps between the representing spectra or spaces.

\begin{description}
\item[Stable.] Stable operations \({F}^*(-) \to {E}^{*+l}(-)\) correspond
to maps of the spectrum \(F\) to \(E\) of degree \(l\), and thus
\(\m{S}^l(F,E) \cong \rc{E}^l(F)\).

\item[Unstable.] Unstable operations \({F}^k(-) \to {E}^l(-)\) correspond
to maps \(\spc{F} \to \spc[l]{E}\) and thus \(\m{U}^l_k(F,E) \cong
{E}^l(\spc{F})\).
\end{description}

\paragraph{Additive Operations.}

Within the family of unstable operations lie the \emph{additive}
operations which we denote by \(\m{A}^l_k(F,E) \subseteq
\m{U}^l_k(F,E)\).  A generic unstable operation need not preserve any
of the structure of \({F}^k(X)\), even that of being an abelian group.
An additive operation is one that does preserve the additive
structure.  Using the fact that the additive structure of \({F}^k(X)\)
comes from the \(H\)\hyp{}space structure of \(\spc{F}\) it is
straightforward to show that within \({E}^l(\spc{F})\) the additive
operations are:
\[
  \ker \left((\mu^* - {p_1}^* - {p_2}^*) \co {E}^l(\spc{F}) \to
  {E}^l(\spc{F} \times \spc{F}) \right)
\]
where \(\mu \co \spc{F} \times \spc{F} \to \spc{F}\) is the \(H\)\hyp{}map
and \(p_1\), \(p_2\) are the projections onto the two factors.  This
is the subspace of \emph{primitives} and is written \(P
{E}^l(\spc{F})\).  Thus \(\m{A}_k^l(F,E) \cong P {E}^l(\spc{F})\).

\paragraph{Co-operations.}

If the spectrum \(F\) and the spaces \((\spc{F})_{k \in \Z}\) have
strong \(E\)\hyp{}duality then the cohomology rings \(\rc{E}^*(F)\)
and \({E}^*(\spc{F})\) are the \(E^*\)\hyp{}duals of the corresponding
homology groups.  Therefore one can analyse the groups of operations
by studying these homology groups.  This is often a Good Thing To Do.
Firstly, the topological issues alluded to in the paragraph on duality
all occur on the cohomology side; homology is discrete.  Secondly, it
is easier to find \emph{explicit} elements in the homology using
push-forwards from key test spaces.

Anything worth studying gets a name, in this case
\emph{co-operations}.  As with operations these come in three
flavours: stable, unstable, and additive.  The stable co-operations
are \(\rc{E}_*(F)\).  The unstable ones are \({E}_*(\spc[*]{F})\).  The
additive co-operations are the \emph{indecomposables} of
\({E}_*(\spc[*]{F})\): for each \(k \in \Z\) we define
\[
  Q{E}_*(\spc{F}) := \coker\big((\mu_* - {p_1}_* - {p_2}_*) \co {E}_*(\spc{F}
  \times \spc{F}) \to {E}_*(\spc{F})\big).
\]
This is a quotient of \({E}_*(\spc{F})\); let \(\tilde{q}_k\) denote the
quotient map.  Assuming sufficient duality the \(E^*\)\hyp{}dual of \(Q
{E}_*(\spc{F})\) is \(P{E}^*(\spc{F})\), which we know to be isomorphic
to \(\m{A}^*_k(F,E)\).

In \cite{jbdjww} the authors regrade the additive co-operations by
defining \(Q(E,F)_*^k := Q{E}_*(\spc{F})\), with the total degree of
\(Q(E,F)_i^k\) being \(k - i\).  The reason for this is that the
algebraic structure of \(Q{E}_*(\spc{F})\) makes more sense with the
new grading.  For this paper there is not much difference between the
two options as we mainly deal with all unstable operations and when we
do explicitly consider additive operations then we are concerned with finding
identities and these, of course, hold whatever the grading scheme in
use.  We choose \(Q(E,F)_*^*\) because \cite{jbdjww} is the main
background for this paper and so we are trying to use their
conventions whenever possible.

The regraded quotient map of degree \(k\) is:
\[
  q_k \co {E}_*(\spc{F}) \to Q(E,F)_*^k.
\]
The stabilisation map \(\sigma_k \co \Sigma^\infty \spc{F} \to F\)
induces the stabilisation map of co-operations:
\[
  {\sigma_k}_* \co E_*(\spc{F}) \to  \rc{E}_*(\spc{F}) \cong
  \rc{E}_*(\Sigma^\infty \spc{F}) \to \rc{E}_{*-k}(F).
\]
This factors through the quotient to additive co-operations.  Thus we
can define the maps:
\begin{align*}
Q{\sigma_k}_* &\co Q {E}_i(\spc{F}) \to \rc{E}_{i-k}(F), \\
Q(\sigma) &\co Q_i^k(E,F) \to \rc{E}_{i-k}(F).
\end{align*}
The former is of degree \(k\), the latter of degree \(0\).

This seems an appropriate place to note that if the spectrum \(F\) and
the spaces \((\spc{F})_{k \in \Z}\) have strong \(E^*\)\hyp{}duality then
the potential \(\lim^1\)\hyp{}problem referred to above disappears: a
stable operation is completely determined by its unstable components.
See \cite[\S 9]{jb4} for more on this issue.

\paragraph{Operations, Maps, and Functionals.}

We therefore have three ways of thinking about operations: as
operations, as maps (or classes), and as functionals on co-operations
(assuming sufficient duality).  We shall distinguish between these
views using fonts and alphabets: roman [italic] for operations, greek
for maps, and gothic for functionals.  We shall attempt to make our
notation as transparent as possible: the stable operation \(r\) will
correspond to the stable map \(\rho\) and to the functional \(\mf{r}\)
on stable co-operations.

In each of the three cases we have a natural restriction map from the
stable to the unstable operations which factors through the additive
ones.  These restriction maps do not correspond exactly: there are
signs to insert at the appropriate junctures.  The full diagram (which
is an expansion of \cite[6.10]{jbdjww}) is:
\begin{equation}
\label{diag:opmapfun}
\begin{array}{c}
  \xymatrix{
    \m{S}^l(F,E) \ar[r] \ar[d]^{\cong} \ar@{}[rd]|{(-1)^{k l}} &
    \m{A}^{k+l}_k(F,E) \ar[r]^{\subseteq} \ar[d]^{\cong} &
    \m{U}^{k+l}_k(F,E) \ar[d]^{\cong} \\
  \rc{E}^l(F) \ar[r]^{{\sigma_k}^*} \ar[d]^{\cong} \ar@{}[rd]|{(-1)^k} &
    P{E}^{k+l}(\spc{F}) \ar[r]^{\subseteq} \ar[d]^{\cong} \ar@{}[rd]|{(-1)^k} &
    {E}^{k+l}(\spc{F}) \ar[d]^{\cong} \\
    D^l\rc{E}_*(F) \ar[r]^{D Q(\sigma)} &
    D^l Q_*^k(E,F) \ar[r]^{D q_k} &
    D^l {E}_*(\spc{F}).
}
\end{array}
\end{equation}

The reasons for the signs in this diagram are quite subtle so it is
worth taking some time to explain them carefully.  This is an
expansion of the \emph{scholium on signs} in \cite[\S6]{jbdjww}.

Let us start with the \((-1)^{k l}\) in the upper left square.  Let
\(r \in \m{S}^l(F,E)\) be a stable operation with restriction \(r_k
\in \m{U}_k^{k+l}(F, E)\).  Consider the following diagram.
\begin{equation}
\label{diag:signs}
\begin{array}{c}
  \xymatrix{
    \rc{F}^0(F) \ar[r]^{{\sigma_k}^*} \ar[d]^{r} \ar@{}[rd]|{(-1)^{k
        l}} &
    \rc{F}^k(\Sigma^\infty \spc{F}) \ar[r]^{\cong} \ar[d]^{r} &
    \rc{F}^k(\spc{F}) \ar[d]^{r} \ar[r]^{\subseteq} &
    {F}^k(\spc{F}) \ar[d]^{r} \\
    \rc{E}^l(F) \ar[r]^{{\sigma_k}^*} &
    \rc{E}^{k+l}(\Sigma^\infty \spc{F}) \ar[r]^{\cong} &
    \rc{E}^{k+l}(\spc{F}) \ar[r]^{\subseteq}&
    E^{k+1}(\spc{F})
}
\end{array}
\end{equation}
As \(r\) is an operation of degree \(l\) and \(\sigma_k\) is a map of
spectra of degree \(k\) we have \({\sigma_k}^* r \alpha = (-1)^{k l} r
{\sigma_k}^* \alpha\).  This accounts for the sign in this diagram.
The other squares commute since the maps involved have degree \(0\).

Let us chase the universal class \(\iota \in \rc{F}^0(F)\) around this
diagram.  The map \(\sigma_k\) was defined so that the image of
\(\iota\) in \({F}^k(\spc{F})\) is \(\iota_k\).  Therefore the image of
\(\iota\) in \({E}^{k+l}(\spc{F})\) when taking the upper route in
diagram~\eqref{diag:signs} is \(r \iota_k\).  As \(r_k\) is the
restriction of \(r\) this is also \(r_k \iota_k\) and so is the class
corresponding to the unstable operation \(r_k\).  Thus the image of
\(\iota\) in \({E}^{k+l}(\spc{F})\) via the upper route in
diagram~\eqref{diag:signs} is the image of \(r\) via the upper route
in diagram~\eqref{diag:opmapfun}.

The lower route starts off with the class \(r \iota \in \rc{E}^l(F)\).
This is the class corresponding to \(r\).  The rest of the lower route
is the map \(\rc{E}^l(F) \to {E}^{k+l}(\spc{F})\) from
diagram~\eqref{diag:opmapfun}.  Hence the image of \(\iota\) via the
lower route in diagram~\eqref{diag:signs} is the image of \(r\) via
the lower of the two routes from \(\m{S}^l(F,E)\) to
\({E}^{k+l}(\spc{F})\) in diagram~\eqref{diag:opmapfun}.

Therefore as we have the sign \((-1)^{k l}\) in
diagram~\eqref{diag:signs} we need it also in
diagram~\eqref{diag:opmapfun}.  As additive operations and classes
are subsets of the unstable ones the sign must go between the stable
and additive lines.

The signs in the lower squares come from a technical subtlety.  Had we
put \(D Q{E}_* (\spc{F})\) in the middle of the lower row there would
have been no signs.  Therefore the signs come from replacing \(D
Q{E}_*(\spc{F})\) by \(D Q_*^k(E,F)\).  Consider the diagram:
\[
  \xymatrix{
    {E}_{k+l}(\spc{F}) \ar[r]^{\tilde{q}_k} \ar[rd]^{q_k} &
    Q{E}_{k+l}(\spc{F}) \ar[r]^{Q {\sigma_k}_*} &
    \rc{E}_l(F) \\
    & Q (E,F)_{k+l}^k \ar[ru]^{Q(\sigma)}  \ar[u]^{\Sigma^{-k}} &
  }
\]
where \(\Sigma^{-k} \co Q (E,F)_*^k \to Q {E}_*(\spc{F})\) is the degree
shift isomorphism and the total upper map is \({\sigma_k}_*\) which
factors as shown.  This diagram commutes.  When we dualise, we find
that:
\begin{align*}
  D \tilde{q}_k &= D(\Sigma^{-k} q_k) = (-1)^k D q_k D (\Sigma^{-k})
  \\
  D Q(\sigma) &= D(Q{\sigma_k}_* \Sigma^{-k}) = (-1)^k D(\Sigma^{-k})
  D Q{\sigma_k}_*
\end{align*}
as in each case both commutants have degree \(k\).  Therefore if we
define the map \(P {E}^{k+l}(\spc{F}) \to D^{k+l}Q (E,F)_*^k\) as:
\[
  P{E}^{k+l}(\spc{F}) \to D^{k+l} Q{E}_*(\spc{F}) \xrightarrow{D
    \Sigma^{-k}} D^l Q (E,F)_*^k
\]
we find that we need to add the signs \((-1)^k\) to each lower square
to make the resultant diagram commute.

\paragraph{Constant Operations and Based Operations.}

There is one particular type of operation that we have to consider, if
only so that we know how to ignore them later.  These are
\emph{constant operations}.  Each \(v \in E^*\) defines an operation
on \(F^*(X)\) by \(x \mapsto v 1_X\), where \(1_X\) is the unit in the
algebra \(E^*(X)\).

Juxtaposed to constant operations are the \emph{based} operations.  An
operation \(r : F^k(-) \to E^l(-)\) is based if it maps zero to zero.
This is, of course, automatic for an additive operation but not for a
general unstable operation.

The reason for mentioning these two types of operation together is
that every (unstable) operation has a decomposition as the sum of a
constant operation and a based operation.  For an operation \(r :
F^*(-) \to E^*(-)\) let \(v_r \in E^* = E^*(\text{pt})\) be the image
of \(0 \in F^* = F^*(\text{pt})\) under \(r\), then
let \(\tilde{r}\) be the based operation given by \(\tilde{r}(x) =
r(x) - v_r 1_X\) for \(x \in F^*(X)\).

The based operations correspond to the classes in
\(\rc{E}^l(\spc{F})\) and thereby to the based maps, \(\lb \spc{F},
\spc[l]{E} \rb_+\).  The based functionals are dual to the reduced
homology groups, \(\rc{E}_*(\spc{F})\).

In each case, the projection from the unbased to the based version is
the obvious one.  Where we have a possibly unbased operation \(r\),
map \(\rho\), or functional \(\mathfrak{r}\) we shall denote the
corresponding based one by \(\tilde{r}\), \(\tilde{\rho}\), or
\(\tilde{\mathfrak{r}}\).

\paragraph{Suspension and Looping.}

There is a method of getting new unstable operations from old.  Given
an unstable operation \(r_k \co {F}^k(-) \to {E}^l(-)\) we can define
another unstable operation \(r_{k-1} \co {F}^{k-1}(-) \to {E}^{l-1}(-)\)
via:
\[
r_{k-1} \co F^{k-1}(X) \to \rc{F}^{k-1}(X) \cong \rc{F}^k(\Sigma X)
\xrightarrow{\tilde{r}_k} \rc{E}^l(\Sigma X) \cong \rc{E}^{l-1}(X)
\subseteq E^{l-1}(X).
\]

The corresponding idea in the world of maps is to use the equivalences
\(\spc[k-1]{E} \simeq \Omega \spc{E}\) and so given a map \(\rho_k \co
\spc{F} \to \spc[l]{E}\) we define \(\rho_{k-1}\) via:
\[
  \rho_{k-1} \co \spc[k-1]{F} \simeq \Omega \spc{F} \xrightarrow{\Omega
    \tilde{\rho}_k} \Omega \spc[l]{E} \simeq \spc[l-1]{E}.
\]

For functionals the push-forward on co-operations defines the
following \emph{suspension} map:
\[
  \Sigma \co E_{l-1}(\spc[k-1]{F}) \to \rc{E}_{l-1}(\spc[k-1]{F})
  \cong \rc{E}_l(\Sigma \spc[k-1]{F})
  \xrightarrow{(-1)^{k-1}{\vartheta_{k-1}}_*} \rc{E}_l(\spc{F})
  \subseteq E_l(\spc{F}).
\]
The sign here is part of the baggage that comes with dealing with
graded and ungraded objects.  Its presence here is a minor nuisance
but its absence would be a minor headache later.  We dualise this map
to one on functionals.

We shall denote this process of getting one operation, map, or
functional from another by \(\Omega\).  Thus, for functionals, \(\Omega
= D \Sigma\).

The diagram relating these maps is:
\[
  \xymatrix{
    \m{U}^l_k(F,E) \ar[r]^{\cong} \ar[d]^{\Omega} \ar@{}[rd]|{1} &
    {E}^l(\spc{F}) \ar[r]^{\cong} \ar[d]^{\Omega} \ar@{}[rd]|{(-1)^k} &
    D^l {E}_*(\spc{F}) \ar[d]^{\Omega} \\
    \m{U}^{l-1}_{k-1}(F,E) \ar[r]^{\cong} &
    {E}^{l-1}(\spc[k-1]{F}) \ar[r]^{\cong} &
    D^{l-1} {E}_*(\spc[k-1]{F})
  }
\]
It is curious that removing the sign from the definition of the
suspension map on functionals does not make this diagram commute
without signs, rather the sign is \(-1\) regardless of the degree.
This is another aspect of the passage from ungraded to graded objects.

We should emphasise that we have defined looping for \emph{unbased}
operations, maps, and functionals.  However, the construction factors
through the projection to the corresponding based objects.

\paragraph{Colimits.}

The spectrum \(F\) is built from the spaces \(\spc{F}\) using the
signed suspension maps \((-1)^k \vartheta_k \co \Sigma \spc{F} \to
\spc[k+1]{F}\).  This expresses \(F\) as equivalent to the colimit of
the sequence \((\Sigma^\infty \spc{F})\) in the category of spectra.
Applying \(E\)\hyp{}homology leads to:
\[
  \rc{E}_*(F) \cong \colim_k \rc{E}_*(\Sigma^\infty \spc[k]{F}) \cong \colim_k
  \rc{E}_*(\spc[k]{F}) \cong \colim_k E_*(\spc{F}).
\]
The last isomorphism is because the suspension map factors through the
projection to reduced homology and so this projection defines an
isomorphism on the colimits.

In particular,
\[
  \rc{E}_l(F) \cong \colim_k \rc{E}_{l+k}(\spc[k]{F}) \cong \colim_k
  E_{l+k}(\spc{F}).
\]
As the suspension map also factors through the quotient to additive
co\hyp{}operations, we can replace \({E}_*(\spc[k]{F})\) by
\(Q(E,F)_*^k\) as appropriate.

\paragraph{Complex Orientation.}

Our cohomology theories are complex orientable so they admit universal
Chern classes.  That is, say for \({F}^*(-)\), there is an element \(\chcl[F]
\in {F}^2(\CP^\infty)\) which restricts to a generator of
\(\rc{F}^*(\CP^1)\) under the canonical inclusion \(\CP^1 \subseteq
\CP^\infty\).  If, identifying once and for all \(\CP^1\) with
\(S^2\), \(\chcl[F]\) restricts to the image of the unit under the natural
isomorphisms \(\rc{F}^{*+2}(S^2) \cong \rc{F}^*(S^0) \cong {F}^*\) then
we say that \(\chcl[F]\) is a \emph{strict universal Chern class}.  Any
universal Chern class can be modified to a strict one so there is
no loss in assuming that all universal Chern classes are strict.

The existence of a universal Chern class implies that
\({F}^*(\CP^\infty) \cong F^* \ldb \chcl[F] \rdb\).
The \(F\)\hyp{}homology of \(\CP^\infty\) is then the free
\(F^*\)\hyp{}module on generators \(\beta^F_i\) of degree \(-2i\)
defined so that \((\chcl[F])^i(\beta^F_j) = \delta^i_j\).

The \(H\)\hyp{}space structure of \(\CP^\infty\) is a map \(\CP^\infty
\times \CP^\infty \to \CP^\infty\).  In cohomology this induces a map:
\[
  F^* \ldb \chcl[F] \rdb \cong {F}^*(\CP^\infty) \to
  {F}^*(\CP^\infty \times \CP^\infty) \cong
   F^* \ldb \chcl[F]_1, \chcl[F]_2 \rdb.
\]
The image of \(\chcl[F]\) under this map is known as the \emph{formal
  group law} of the cohomology theory \({F}^*(-)\).  We shall write
this formal power series as
\[
  \chcl[F]_1 \plusf \chcl[F]_2
\]
(the ``\(F\)'' is to indicate the
cohomology theory).

In certain circumstances it is possible to substitute elements of an
\(F^*\)\hyp{}algebra into the formal power series that this represents.
(The only difficulty here is with convergence of the resulting sum; so
it works, for example, on nilpotent elements and it works if the
algebra is complete with respect to some filtration and successive
powers of the elements that one is substituting in lie further and
further down in the filtration.)  The properties of the formal group
law imply that, when this is possible, the resulting operation is
associative, commutative, unital, and has inverses \hyp{} hence the name
``formal group law''.  We shall denote iterations of this process with
the adorned summation notation:
\[
  \sideset{}{^F}\sum
\]

We shall need one more fact about the structure of the formal group
law as a power series.  It follows from the basic properties of formal
group laws that
there are identities:
\begin{equation}
\label{eq:fglred}
  \chcl[F]_1 \plusf \chcl[F]_2 = \chcl[F]_1 + \chcl[F]_2 R_1(\chcl[F]_1, \chcl[F]_2) = \chcl[F]_2 + \chcl[F]_1
  R_2(\chcl[F]_1, \chcl[F]_2)
\end{equation}
for some formal power series \(R_1(\chcl[F]_1, \chcl[F]_2), R_2(\chcl[F]_1, \chcl[F]_2)\).

A particular case where substitution is allowed is the element
\(\chcl[F]\) of \(F^* \ldb \chcl[F] \rdb\).  Substituting this into
both variables we define
\[
  \pseries{2}{F}(\chcl[F]) = \chcl[F] \plusf \chcl[F] \in F^*\ldb \chcl[F] \rdb.
\]
It is straightforward to see that the resulting formal power series
has leading term \(2 \chcl[F]\) and so can be again substituted in to
the formal group law.  Iterating this procedure, we define
\(\pseries{n}{F}(\chcl[F]) := \chcl[F] \plusf
\pseries{n-1}{F}(\chcl[F])\).  This formal power series is called the
\emph{\(n\)\hyp{}series} of \({F}^*(-)\).

There is an alternative derivation of these formal power series.  The
\(H\)\hyp{}space structure on \(\CP^\infty\) defines an \(n\)th power
map \(\CP^\infty \to \CP^\infty\).  Using the isomorphism
\({F}^*(\CP^\infty) \cong F^* \ldb \chcl[F] \rdb\), the image of
\(\chcl[F]\) under the pull-back via this map is a formal power series
in \(\chcl[F]\) and it is not hard to see that it is
\(\pseries{n}{F}(\chcl[F])\).

A particularly important case of this is the \emph{\(p\)\hyp{}series} for
\(p\) a prime.  This is of the form:
\[
  \pseries{p}{F}(\chcl[F]) = p \chcl[F] + \sum_{j \ge 1} g^F_j (\chcl[F])^{j+1}
\]
for some \(g^F_j \in {F}^{-2 j}\).  The reduction of this modulo \(p\)
has the form:
\[
  \pseries{p}{F}(\chcl[F]) \equiv \sideset{}{^F}\sum_{i \ge 1} v^F_i
  (\chcl[F])^{p^i} \mod p
\]
for some \(v^F_i \in {F}^{-2(p^i - 1)}\).  Note the adorned summation
sign.

The Chern class for \({F}^*(-)\) is represented by a map \(\chcl[F]
\co \CP^\infty \to \spc[2]{F}\).  Applying \(E\)\hyp{}homology leads
to a push-forward \(\chcl[F]_{\:*} \co {E}_*(\CP^\infty) \to
{E}_*(\spc[2]{F})\).  As \({E}^*(-)\) is itself complex
orientable the former is the free \(E^*\)\hyp{}module on generators
\(\beta^E_i\).  Let \(b_i = \chcl[F]_{\:*} \beta^E_i\) and define:
\[
  b(s) = \sum_{i \ge 0} b_i s^i \in {E}_*(\spc[*]{F}) \ldb s \rdb.
\]

We shall use the same notation, i.e.~\(b_i\), for the images of the
\(b_i\) in the additive and stable co-operations.

\subsection{Algebraic Structure}

The various groups of operations and co-operations have considerable
algebraic structure.  The full list is long so we shall only describe
what we need.  For all the gory details, see \cite{jb4} and
\cite{jbdjww}.

The main structures that we shall use are the multiplicative and
bimodule structures on the sets of co-operations and the bimodule
structure on the sets of operations.  This is further complicated by
the fact that there are two multiplications on the unstable
co-operations.

Once we have introduced these algebraic structures we shall consider
how some of the data we have already seen behaves algebraically.

\paragraph{Co-operation Multiplications.}

The more important \hyp{} for our purposes \hyp{} multiplication is defined
using the maps on the spaces \(\spc{F}\) and spectrum \(F\) which
represent the multiplication in \({F}^*(-)\).  That is, the
map \(\phi_{l,k} \co \spc[l]{F} \times \spc{F} \to \spc[l + k]{F}\)
defines a push-forward:
\[
  {E}_*(\spc[l]{F}) \times {E}_*(\spc{F}) \to {E}_*(\spc[l]{F} \times
  \spc{F}) \xrightarrow{{\phi_{l,k}}_*} {E}_*(\spc[l+k]{F}).
\]
As \(\phi_{l,k}\) is a component of an infinite loop map we also get
multiplications on the additive and stable sets of co-operations which
all correspond under the maps from unstable co-operations to additive
and to stable.  For unstable co-operations we shall write this
multiplication as \((a,b) \mapsto a \circ b\).  For the others we
shall just use the abutment notation.  Note that as the quotient from
unstable to additive co-operations has a non-trivial degree, the
correct formula on a product is:
\[
  q_{i+j}(a \circ b) = (-1)^{j \abs{a}}q_i(a) q_j(b)
\]
for \(a \in {E}_*(\spc[i]{F})\) and \(b \in {E}_*(\spc[j]{F})\).

For additive and stable co-operations these multiplications are graded
commutative (taking the total degree in the regraded additive realm).
For unstable co-operations this is still true but the issue is
somewhat complicated by the fact that the set of unstable
co-operations, \({E}_*(\spc[*]{F})\), has two indices which are used in
different ways: the first is a genuine grading whereas the second is
really only a labelling.  However this multiplication does use this
second index.  To describe exactly how, we would need to introduce yet
more of the structure and it turns out that, for our purposes, this is
unnecessary since any element with both indices even commutes with
everything.  On the few occasions where we need to consider other
elements we shall give the explicit commutation formula.

In light of this confusion, we add that when we speak of the degree of
an element in \({E}_*(\spc[*]{F})\) we shall be using the first index
only.

The set of unstable co-operations has another multiplication which
comes from the \(H\)\hyp{}map \(\spc{F} \times \spc{F} \to
\spc{F}\).  This is graded commutative with the ``honest'' grading.
Note that this product only makes sense for elements which have the
same \emph{second} index.  We shall write this multiplication as
\((a,b) \mapsto a * b\).

The interaction of the two multiplications is controlled by a
coproduct, \(\psi\), which is induced by the diagonal map \(\spc{F}
\to \spc{F} \times \spc{F}\).  That is, if \(\psi(c) = \sum_i c_i'
\otimes c_i''\) then:
\[
  (a * b) \circ c = \sum_i (-1)^{\abs{b} \abs{\right.c_i'\left.}} (a
  \circ c_i') * (b \circ c_i'').
\]
This is the only place where we use this coproduct.

The reason that the \(*\)\hyp{}product does not appear in the additive or
stable realms is that it is what is being quotiented out when passing
to the additive co-operations.  Specifically, the quotient on a
\(*\)\hyp{}product is:
\[
q_k(a * b) = \epsilon_k(a) b + (-1)^{\abs{a} \abs{b}} \epsilon_k(b)a,
\]
where \(\epsilon_k\) is the appropriate augmentation.

\paragraph{Bimodule Structure.}

The various groups of operations from \(F\)\hyp{}cohomology to
\(E\)\hyp{}cohomology have the structure of \(({E}^*\!\! -\!\!
{F}^*)\)\hyp{}bimodules.  
The left \(E^*\)\hyp{}action is:
\[
  (v \cdot r)(\alpha) = v r(\alpha)
\]
whilst the right \({F}^*\)\hyp{}action is:
\[
  (r \cdot v)(\alpha) = r(v \alpha).
\]

In terms of maps these actions are given by composition with certain
maps of the representing spaces.  For \(v \in {E}^l = {E}^l(\text{pt})\)
we define \(\xi v \co \spc{E} \to \spc[k+l]{E}\) by:
\[
  \spc{E} \cong \text{pt} \times \spc{E} \xrightarrow{v \times 1}
  \spc[l]{E} \times \spc{E} \xrightarrow{\phi_{l,k}} \spc[k+l]{E}.
\]
In the stable case we use the smash product and view \(v\) as an
element of \(\rc{E}^l(S)\) (we could have used the smash product in
the unstable case as well since the multiplication factors through the
smash product).  Using these maps we define the left action of \(E^*\)
and right action of \({F}^*\) by appropriate composition:
\[
  v \cdot \rho = (\xi v) \rho, \qquad \rho \cdot v = \rho (\xi v).
\]
The left action of \(E^*\) agrees with the obvious action on
\({E}^*(\spc{F})\).

For co-operations we have an obvious left action of \(E^*\) as the
coefficient ring.  The right action of \({F}^*\) is given by 
push-forwards:
\[
  (\xi v)_* \co {E}_*(\spc{F}) \to {E}_*(\spc[k+l]{F}).
\]
Unpacking the construction of \(\xi v\), and using the definition of
the \(\circ\)\hyp{}\hspace{0pt}multiplication, we see that there is an
element \(\lb v \rb \in {E}_0(\spc[l]{F})\) such that the right action
of \(v\) on \({E}_*(\spc{F})\) is: \(c \mapsto c \circ \lb v \rb\).  The
element \(\lb v \rb\) is the image of \(1\) under the map \(v_* \co {E}^*
= {E}_*(\text{pt}) \to {E}_*(\spc[l]{F})\).  There are corresponding
actions in the additive and stable realms since the map \(\xi v\) is a
component of an infinite loop map.

Diagram~\eqref{diag:opmapfun} is then a diagram of
\(({E}^*\!\!-\!\!{F}^*)\)\hyp{}bimodules.

\paragraph{Algebraic Suspension.}

The suspension map on functionals has a particularly pleasant
structure.  The suspension isomorphism \(\rc{E}_0(S^0) \cong
\rc{E}_1(S^1)\) defines a canonical element \(u_1 \in \rc{E}_1(S^1)\)
as the image of the unit.   This element determines the
suspension isomorphism as follows.  The \(E^*\)\hyp{}module
\(\rc{E}_*(S^1)\) is free of rank one generated by \(u_1\) so we have the
following isomorphisms:
\begin{equation}
\label{eq:uone}
\rc{E}_k(\Sigma X) = \rc{E}_k(S^1 \smsh X) \cong \left(\rc{E}_*(S^1)
  \otimes_{{E}^*} \rc{E}_*(X) \right)_k \cong \rc{E}_{k-1}(X)
\end{equation}
where the final map is \(u_1 \otimes c \mapsto c\).

From equation~\eqref{eq:unit}, the map \(\vartheta_{k-1} \co \Sigma
\spc[k-1]{F} \to \spc{F}\) factors as \(\phi_{1, k-1}(\eta_1 \smsh
1)\).  Thus the following diagram commutes:
\[
  \xymatrix@C=13pt@R=24pt{
    \rc{E}_{l-1}(\spc[k-1]{F}) \ar [r]^-{\cong} &
    \Big(\rule{0pt}{14.8pt} \rc{E}_*(S^1) \otimes_{{E}^*}
    \rc{E}_*(\spc[k-1]{F})\Big)_{l} \ar[r]^-{\cong}
    \ar[d]^{{\eta_1}_* \otimes 1} &
    \rc{E}_{l}(S^1 \smsh \spc[k-1]{F}) \ar[r]^-{\cong}
    \ar[d]^{(\eta_1 \smsh 1)_*} &
    \rc{E}_{l}(\Sigma \spc[k-1]{F}) \ar[dd]^{{\vartheta_{k-1}}_*} \\
    &
    \Big(\rule{0pt}{14.8pt} \rc{E}_*(\spc[1]{F}) \otimes_{{E}^*} \rc{E}_*(\spc[k-1]{F})
    \Big)_{l} \ar[r] &
    \rc{E}_{l}(\spc[1]{F} \smsh \spc[k-1]{F}) \ar[d]^{{\phi_{1,k-1}}_*}
    \\
    &&
    \rc{E}_l(\spc{F}) \ar[r]^{=} &
    \rc{E}_l(\spc{F}) \\
}
\]

The upper route is, up to sign, the suspension map.  Thus
from the lower route, we can see that this map is:
\[
  c \mapsto (-1)^{k-1} e \circ c
\]
where \(e = {\eta_1}_* u_1 \in \rc{E}_1(\spc[1]{F})\).  We shall use
the same notation for the image of \(e\) in \(Q(E,F)_1^1\).  In the
stable realm it maps to the identity (the maps which define stable
co-operations as the colimit of unstable are, up to sign,
\(\circ\)\hyp{}multiplication by \(e\)).

The commutation law for, coproduct of, and augmentation of the element
\(e\) are:
\begin{align*}
  a \circ e &= (-1)^{j + k} e \circ a, \qquad a \in {E}_j(\spc{F}); \\
  \psi(e) &= e \otimes 1_1 + 1_1 \otimes e; \\
  \epsilon_1(e) &= 0.
\end{align*}

\paragraph{Algebraic Chern Class.}

Returning to the series \(\sum b_i s^i\), the first two terms are
readily identifiable in terms of the algebraic structure.  The first,
\(b_0\), is \(1_2\), the \(*\)\hyp{}unit in \({E}_0(\spc[2]{F})\).  The
second, as our Chern classes were strict, is \(- e^{\circ 2}\).  These
quotient to the (regraded) additives as follows:
\begin{align*}
  q_2(b_0) &= q_2(1_2) =   0 \\
  q_2(b_1) &= q_2(-e \circ e) = q_1(e) q_1(e) = e^2.
\end{align*}

The \(b_i\) \(\circ\)\hyp{}commute with everything as they lie in
\({E}_{2i}(\spc[2]{F})\).  Their coproducts and augmentations are:
\begin{align*}
\psi(b_k) &= \sum_{i+j = k} b_i \otimes b_j; \\
\epsilon_2(b_k) &= \begin{cases}
0 & \text{if } k > 0, \\
1 & \text{if } k = 0.
\end{cases}
\end{align*}

\subsection{Morava K\hyp{}Theory}

The Morava K\hyp{}theories will be our main examples of target
theories.  These are a family of multiplicative generalised cohomology
theories indexed by primes and non-negative integers.  There are some
peculiarities corresponding to prime \(2\) which we wish to avoid so
we fix an odd prime, \(p\).  For any prime the theory corresponding to
zero is ordinary rational cohomology so any interesting behaviour
peculiar to the Moravian theories would be expected to rear its head
for strictly positive integers, and this is true for the phenomenon we
have observed, hence we choose \(n \ge 1\).  Thus we have fixed our
attention on \(K(n)^*(-)\), the \(n\)th Morava K\hyp{}theory at the
prime \(p\), for \(n > 0\) and \(p\) odd.  (The prime is not explicit
in the notation as it is quite unusual to vary it in the course of a
discussion whereas it is sometimes fruitful to consider
different values of \(n\).)

The coefficient ring of \(\kn^*(-)\) is
\[
  \kn^* = \F_p\lb {v_n}, v_n^{-1} \rb
\]
where \(\abs{v_n} = -2(p^n-1)\).  This is a graded field and hence all
modules over this ring are free.  Two consequences of this are that
\(\kn^*(-)\) has a K\"unneth formula and has strong duality.

The \(p\)\hyp{}series for \(\kn^*(-)\) is
\[
  \lb p \rb_{\kn}(s) = v_n s^{p^n}.
\]

\section{Analysing the \protect\(p\protect\)\hyp{}Series}
\label{sec:pseries}

In this section we analyse what information can be gleaned from the
\(p\)\hyp{}series of the two cohomology theories under consideration.
From now on we assume that \({E}^*(-)\) and \({F}^*(-)\) satisfy the
conditions of theorem~\ref{th:main}.  That is, they are multiplicative
graded cohomology theories which are commutative and complex
orientable and the following conditions hold.
\begin{enumerate}
\item The coefficient ring of \({E}^*(-)\) has characteristic \(p\).

\item The formal group law of \({E}^*(-)\) has finite height, say \(n\).

\item The coefficient of the first term in the \(p\)\hyp{}series for
\({E}^*(-)\) is invertible.

\item The various groups of operations from \({F}^*(-)\) to \({E}^*(-)\)
are dual over the coefficient ring of \({E}^*(-)\) to the corresponding
groups of co-operations.
\end{enumerate}

The main tool in our analysis is a result from
\cite{drww}.
\begin{theorem}[Ravenel-Wilson]
The following identity holds in \({E}_*(\spc[*]{F})\ldb s \rdb\):
\[
  b(\lb p \rb_E(s)) = \lb p \rb_F(b(s)),
\]
where, in expanding out the right-hand side, the coefficients
\(g_j^F\) of the \(p\)\hyp{}series for \({F}^*(-)\) act via the right action
of \({F}^*\) on \({E}_*(\spc[*]{F})\).
\end{theorem}
Recall that \(b(s) = \sum_{i \ge 0} b_i s^i\).

To unpack this we use the fact that the maps which represent the
addition and multiplication in \({F}^*(-)\) defined the \(*\)\hyp{} and
\(\circ\)\hyp{}multiplications on \({E}_*(\spc[*]{F})\).  Therefore, when
expanding the right-hand side, we need to translate addition to
\(*\)\hyp{}multiplication and multiplication to \(\circ\)\hyp{}multiplication.
This leads to:
\begin{equation}
\label{eq:pseries}
  b\left( p s + \sum_{i > 0} g_i^E s^{i+1} \right) = b(s)^{* p} *
  \bstar_{i > 0} \left( b(s)^{\circ i + 1} \circ \lb g_i^F \rb\right).
\end{equation}

\subsection{Additive Co-operations}

Equation~\eqref{eq:pseries} looks horrendous but simplifies
considerably when we quotient to the additive co-operations.
Throughout this section we shall be working in the additive realm;
that is, with \(Q (E,F)^*_*\) and formal power series over this.

As \(\epsilon b(s) = 1\), we find that in \(Q (E,F)^*_* \ldb s \rdb\):
\[
  b \left( p s + \sum_{i > 0} g_i^E s^{i+1} \right) = p b(s) + \sum_{i
    > 0} b(s)^{i+1} \lb g_i^F \rb.
\]

In the additive realm it is a tautology that the left and right
\Z{}\hyp{}actions agree.  As \({E}_*(\spc[*]{F})\) is an
\(E^*\)\hyp{}module it has characteristic \(p\) and thus we may
replace both sides by their reductions modulo \(p\).  That is:
\begin{equation}
\label{eq:equatep}
  b \left( \sideset{}{^E}\sum_{i > 0} v_i^E s^{p^i} \right) =
  \sideset{}{^F}\sum_{i > 0} b(s)^{p^i} \lb v_i^F \rb.
\end{equation}

From this equation we shall deduce the following result.

\begin{proposition}
\label{prop:addloop}
For \(n \in \N\), let \(\pi_n = \frac{p^n-1}{p-1}\).  Then in \(Q
(E,F)_*^*\):
\begin{equation}
\label{eq:addloop}
v_n^E {b_1}^{\pi_n} = {b_1}^{\pi_{n+1} - 1} \lb v_n^F \rb.
\end{equation}
\end{proposition}

\begin{proof}
\renewcommand{\qed}{}
Our strategy for proving \eqref{eq:addloop} is to equate powers of
\(s\) in \eqref{eq:equatep} and read off certain identities.  To
begin, we examine the left-hand side of \eqref{eq:equatep} to find its
leading term.  The left-hand side is of the form \(b(r(s))\) where
\(r(s) = \sum^E_{i > 0} v_i^E s^{p^i}\).  As \(b(s)\) has
leading term \(b_1 s\), the leading term of \(b(r(s))\) is the product
of \(b_1\) and the leading term of \(r(s)\).

To find this leading term we use the formula in \eqref{eq:fglred}.
Let \(r_i(s)\), \(i \in \{1,2\}\), be formal power series in \(s\)
with leading terms \(a_i s^{l_i}\) and suppose that \(1 \le l_1 <
l_2\).  Then \eqref{eq:fglred} shows that
\begin{align*}
  r_1(s) \plusf[E] r_2(s) &= r_1(s) \mod s^{l_2}\\
  &= a_1 s^{l_1} \mod s^{l_1
    + 1}.
\end{align*}
As \(r(s)\) is a summation (using the formal group law of \({E}^*(-)\))
of monomials of strictly increasing degree, the above shows that its
leading term will be the first non-zero monomial.  Our assumptions on
the cohomology theory \({E}^*(-)\) imply that this is \(v_n^E s^{p^n}\).
Hence the leading term of the left-hand side of \eqref{eq:equatep} is
\(v_n^E b_1 s^{p^n}\).

Now let us consider the right-hand side of \eqref{eq:equatep}.  As
\(b(s)\) has leading term \(b_1 s\), \(b(s)^{p^i}\) has leading term
\({b_1}^{p^i} s^{p^i}\).  Therefore the above argument shows that the
leading term of the right-hand side of \eqref{eq:equatep} is \({b_1}^p
s^p \lb v_1^F \rb\).  If \(n = 1\), equating coefficients of \(s^p\)
yields the identity:
\[
  v_n^E b_1 = {b_1}^p \lb v_1^F \rb.
\]
This is precisely \eqref{eq:addloop} with \(n = 1\) since \(\pi_1 =
1\) and \(\pi_2 = p+1\).

If \(n \ne 1\), equating coefficients yields:
\[
  0 = {b_1}^p \lb v_1^F \rb.
\]

This provides the start of a recursion procedure which will lead to
our desired result.  Assume that \(n > 1\) and that for some \(m \in
\N\) with \(1 < m < n\) we have
\[
  0 = {b_1}^{\pi_{j+1} - 1} \lb v_j^F \rb
\]
for all \(j \in \N\) such that \(1 \le j < m\).  As \(\pi_2 = p + 1\) we
have shown above that this holds for \(m = 2\).

We multiply \eqref{eq:equatep} through by \({b_1}^{\pi_{m} - 1}\).  The
leading term of the left-hand side of the resulting equation is simply
\(v_n^E {b_1}^{\pi_{m}} s^{p^n}\).  Let us examine the right-hand side.
We apply another recursion argument.  Suppose that for some \(1 \le j
< m\) we have
\[
  {b_1}^{\pi_{m} - 1} \sideset{}{^F}\sum_{i > 0} b(s)^{p^i} \lb v_i^F
  \rb =
  {b_1}^{\pi_{m} - 1} \sideset{}{^F}\sum_{i \ge j} b(s)^{p^i} \lb v_i^F
  \rb.
\]
Note that when \(j = 1\) this is a tautology.  Using \eqref{eq:fglred}
we expand this out:
\begin{align*}
  {b_1}^{\pi_{m} - 1} \sideset{}{^F}\sum_{i > 0} b(s)^{p^i} \lb v_i^F
  \rb &=
  {b_1}^{\pi_{m} - 1} \sideset{}{^F}\sum_{i \ge j} b(s)^{p^i} \lb v_i^F
  \rb \\
  &= {b_1}^{\pi_{m} - 1} \left( b(s)^{p^j} \lb v_j^F \rb \plusf
  \sideset{}{^F}\sum_{i > j} b(s)^{p^i} \lb v_i^F \rb \right) \\
  &= {b_1}^{\pi_{m} - 1} \left( \Bigg(\sideset{}{^F}\sum_{i > j} b(s)^{p^i}
  \lb v_i^F \rb \Bigg) + (b(s)^{p^j} \lb v_j^F \rb) R(s) \right) \\
\intertext{\hfill for some formal power series \(R(s) \in Q
    (E,F)^*_* \ldb s \rdb\)}
  &= {b_1}^{\pi_{m} - 1} \Bigg( \sideset{}{^F}\sum_{i > j} b(s)^{p^i}
  \lb v_i^F \rb \Bigg) + b(s)^{p^j} {b_1}^{\pi_{m} - 1} \lb v_j^F \rb
  R(s) \\
  &= {b_1}^{\pi_{m} - 1} \sideset{}{^F}\sum_{i \ge j+1} b(s)^{p^i} \lb
  v_i^F \rb.
\end{align*}
The last line follows since \(j < m\) and, by assumption,
\({b_1}^{\pi_{j+1} - 1} \lb v_j^F \rb = 0\).

The last time we can apply our recursion is when \(j = m - 1\).  This
leaves us with the equation
\[
  {b_1}^{\pi_{m} - 1} \sideset{}{^F}\sum_{i > 0} b(s)^{p^i} \lb v_i^F
  \rb =
  {b_1}^{\pi_{m} - 1} \sideset{}{^F}\sum_{i \ge m} b(s)^{p^i} \lb
  v_i^F \rb.
\]
By a now-familiar argument, the leading term of the right-hand side of
this is
\[
  {b_1}^{\pi_{m} - 1} {b_1}^{p^{m}} \lb v_{m}^F \rb s^{p^{m}} =
  {b_1}^{\pi_{m+1} - 1} \lb v_{m}^F \rb s^{p^{m}}.
\]

We now equate coefficients in the modified version of
\eqref{eq:equatep}, i.e.~after multiplying both sides by
\({b_1}^{\pi_{m} - 1}\).  If \(m < n\) we see that
\[
0 = {b_1}^{\pi_{m+1} - 1} \lb v_{m}^F \rb
\]
whence we can continue our recursion.

This recursive argument stops when \(m = n\) for then we no longer
have zero on the left-hand side.  Equating coefficients at this point
yields the desired equation

\hspace*{\fill}
\(\displaystyle
 v_n^E {b_1}^{\pi_n} = {b_1}^{\pi_{n+1} - 1} \lb v_n^F \rb.\)
\hspace*{\fill}\makebox[0pt][r]{\qedsymbol}
\end{proof}

Since \(\pi_{n+1} = p^n + \pi_n\) we have shown that if \(h = \pi_n\)
then the following holds in \(Q (E,F)^*_*\):
\[
  v_n^E {b_1}^h = {b_1}^{p^n + h - 1} \lb v_n^F \rb.
\]
It is entirely possible that this will hold for some smaller value of
\(h\) and the minimum such value is an interesting invariant of the
cohomology theory \({F}^*(-)\).  In light of the fact that \(b_1 = e^2\)
we get slightly finer control if we consider this as an identity about
\(e\) rather than \(b_1\).

\begin{defn}
\label{def:addheight}
Let \({E}^*(-)\) and \({F}^*(-)\) be complex orientable, graded,
commutative, multiplicative cohomology theories.  Suppose that the
coefficient ring, \(E^*\), has characteristic \(p\) and that the
formal group law for \({E}^*(-)\) has finite height, say \(n\).  Let
\(v_n^E \in {E}^{-2(p^n-1)}\) and \(v_n^F \in {F}^{-2(p^n-1)}\) be the
coefficients of \(s^{p^n}\) that appear in the formal sum giving the
mod \(p\) reduction of the \(p\)\hyp{}series of the formal group laws of
\({E}^*(-)\) and \({F}^*(-)\) respectively.

Define the \emph{\(E\)\hyp{}additive loop height} of \({F}^*(-)\) to be the
least positive integer \(h\) for which the identity
\[
  v_n^E e^h = e^{2(p^n-1) + h} \lb v_n^F \rb
\]
holds in \(Q (E,F)^*_*\).

In the case that \({F}^*(-) = {E}^*(-)\) we shall refer to this as the
\emph{self additive loop height} of \({E}^*(-)\).
\end{defn}

\begin{examples}
\begin{enumerate}
\item By proposition~\ref{prop:addloop} the maximum possible
\(E\)\hyp{}additive loop height is \(2 \frac{p^n-1}{p-1}\) where \(n\) is
the height of the formal group law of \({E}^*(-)\).  The distinct lack of
any relations in \(\kn_*(\spc[*]{\bp})\), as demonstrated in
\cite{drww}, allows one to conclude that the \(\kn\)\hyp{}additive loop
height for \(\bp\) is \(2 \frac{p^n-1}{p-1}\).  Thus the bound given
in proposition~\ref{prop:addloop} is the best possible.

\item On the other hand, \cite[Proposition 1.1(j)]{ww} implies the
self additive loop height of \(\kn^*(-)\) is \(1\).
\end{enumerate}
\end{examples}

\subsection{Unstable Co-operations}

The analysis of the \(p\)\hyp{}series in the unstable realm follows from
that in the additive context due to a very useful trick:
\(\circ\)\hyp{}multiplication by \(e\) factors through additive
co-operations.  That is, if we have unstable co-operations \(a, c\)
such that \(q_k(a) = q_k(c)\) then \(e \circ a = e \circ c\).  Thus we
can ignore equation~\eqref{eq:pseries} and apply the additive results
to the unstable situation.

We have an unstable version of definition~\ref{def:addheight}:
\begin{defn}
\label{def:unstabheight}
Let \({E}^*(-)\) and \({F}^*(-)\) be complex orientable, graded,
commutative, multiplicative cohomology theories.  Suppose that the
coefficient ring, \(E^*\), has characteristic \(p\) and that the
formal group law for \({E}^*(-)\) has finite height, say \(n\).  Let
\(v_n^E \in {E}^{-2(p^n-1)}\) and \(v_n^F \in {F}^{-2(p^n-1)}\) be the
coefficients of \(s^{p^n}\) that appear in the formal sum giving the
mod \(p\) reduction of the \(p\)\hyp{}series of the formal group laws of
\({E}^*(-)\) and \({F}^*(-)\) respectively.

Define the \emph{\(E\)\hyp{}unstable loop height} of \({F}^*(-)\) to be the
least positive integer \(h\) for which the identity
\[
  v_n^E e^{\circ h} = e^{ \circ (2(p^n-1) + h)} \circ \lb v_n^F \rb
\]
holds in \({E}_*(\spc[*]{F})\).

In the case that \({F}^*(-) = {E}^*(-)\) we shall refer to this as the
\emph{self unstable loop height} of \({E}^*(-)\).
\end{defn}

The argument above produces:
\begin{lemma}
The \(E\)\hyp{}unstable loop height of \({F}^*(-)\) is at least the
\(E\)\hyp{}additive loop height and at most one more.  In particular, \(2
\frac{p^n-1}{p-1} + 1\) is an upper bound. \noproof
\end{lemma}

Careful examination of \cite[Proposition 1.1(j)]{ww} reveals that the
self unstable loop height of \(\kn^*(-)\) is \(1\).

\section{Splitting Co-operations}
\label{sec:splitting}

In this section we use the results of the previous one to define how
to construct a stable operation from an unstable one.  Our strategy
will be to use the formula from proposition~\ref{prop:addloop}, and
its unstable version, to define idempotents in the co-operation
algebras which will split the co-operations.

\subsection{Idempotents}

\begin{defn}
Let \(s \in {E}_0(\spc[0]{F})\) denote the unstable co-operation:
\[
  s := (v_n^E)^{-1} e^{\circ 2(p^n - 1)} \circ \lb v_n^F \rb.
\]
\end{defn}

Recall that one of the conditions on the cohomology theory \({E}^*(-)\)
is that the element \(v_n^E \in {E}^*\) is invertible, hence \(s\) is
well-defined.

\begin{proposition}
\label{prop:idemp}
Let \(h\) be the \(E\)\hyp{}unstable loop height of \({F}^*(-)\).  The
co-operation \(s\) has the following properties:
\begin{enumerate}
\item \(s \circ s = s\); that is, \(s\) is an idempotent for
the \(\circ\)\hyp{}multiplication.

\item \(e \circ s = s \circ e\).

\item \(e^{\circ h} \circ s = e^{\circ h} = s \circ e^{\circ h}\).

\item There is some \(s'\) such that \(e^{\circ h} \circ s' =
s\).
\end{enumerate}
\end{proposition}

\begin{proof}
\renewcommand{\qed}{}
\begin{enumerate}
\item As \(h\) is the \(E\)\hyp{}unstable loop height of \({F}^*(-)\) we
have the identity:
\[
  v_n^E e^{\circ h} = e^{\circ (2(p^n-1) + h)} \circ \lb v_n^F \rb
\]
which rearranges to:
\[
  e^{\circ h} = (v_n^E)^{-1} e^{\circ (2(p^n-1) + h)} \circ \lb v_n^F
  \rb.
\]
Now \(h \le 2\frac{p^n - 1}{p-1} + 1\).  As \(p\) is odd, it is at
least \(3\) and so \(h\) is strictly less than \(2(p^n-1)\).  Hence:
\[
  e^{\circ 2(p^n-1)} = (v_n^E)^{-1} e^{\circ 2(p^n - 1)} \circ
  e^{\circ 2(p^n - 1)} \circ \lb v_n^F \rb,
\]
which leads to:
\[
  (v_n^E)^{-1} e^{\circ 2(p^n-1)} \circ \lb v_n^F \rb = (v_n^E)^{-1}
  e^{\circ 2(p^n-1)} \circ \lb v_n^F \rb \circ (v_n^E)^{-1} e^{\circ
    2(p^n-1)} \circ \lb v_n^F \rb.
\]
This is another way of saying that \(s \circ s = s\).

\item As \(s\) has both indices zero it \(\circ\)\hyp{}commutes with
everything.

\item From
\[
  e^{\circ h} = (v_n^E)^{-1} e^{\circ (2(p^n-1) + h)} \circ \lb v_n^F
  \rb
\]
we deduce that
\[
  e^{\circ h} = e^{\circ h} \circ s.
\]
Then \(s \circ e^{\circ h} = e^{\circ h}\) as \(s\)
\(\circ\)\hyp{}commutes with everything.

\item As \(h < 2(p^n-1)\) the element \(s' = (v_n^E)^{-1} e^{\circ
  (2(p^n-1) - h)} \circ \lb v_n^F \rb\) is well-defined.  It clearly
has the desired property. \hspace*{\fill}\qedsymbol
\end{enumerate}
\end{proof}

\begin{corollary}
There is a split short exact sequence of graded algebras (using the
\(\circ\)\hyp{}multiplication):
\[
  0 \to s {E}_*(\spc[*]{F}) \to {E}_*(\spc[*]{F}) \to {E}_*(\spc[*]{F})/s
  {E}_*(\spc[*]{F}) \to 0.
\]
The first splitting map is \(\circ\)\hyp{}multiplication by \(s\).  The second
identifies the quotient algebra with the ideal generated by \((1 -
s)\).

Let \(h\) be the \(E\)\hyp{}unstable loop height of \({F}^*(-)\).  The map
\(\Sigma^h\) is an isomorphism on the ideal generated by \(s\) and is
null on the ideal generated by \((1 - s)\).
\end{corollary}

\begin{proof}
This is essentially a rephrasing of proposition~\ref{prop:idemp} in
terms of maps rather than elements.  Define a map:
\[
  S \co {E}_*(\spc[*]{F}) \to {E}_*(\spc[*]{F}), \qquad S(c) = s \circ c.
\]
As \(s\) is an idempotent, \(S\) is a projection and an algebra
map.  The splitting follows by basic algebra.

Up to sign, the map \(\Sigma^h\) is multiplication by \(e^{\circ h}\).
Let \(S'\) be the operation of \(\circ\)\hyp{}multiplication by \(s'\).
The two latter properties of \(s\) show that: \(\Sigma^h S =
\Sigma^h = S \Sigma^h\) and \(S' \Sigma^h = S = \Sigma^h S'\).  From
these we readily see that \(\im \Sigma^h = \im S\) and \(\ker \Sigma^h
= \ker S\).  Thus \(\Sigma^h\) restricts to an isomorphism on the
image of \(S\) and is null on the kernel.
\end{proof}

\begin{corollary}
All of the above quotients to the additive realm. \noproof
\end{corollary}

\subsection{Colimits}

We label the various maps in the split short exact sequence as
follows:
\[
  0 \to s {E}_*(\spc[*]{F}) \xrightarrow{\iota_S} {E}_*(\spc[*]{F})
  \xrightarrow{\pi_{\hat{S}}} {E}_*(\spc[*]{F})/s {E}_*(\spc[*]{F}) \to 0
\]
with splitting maps \(\pi_S\) and \(\iota_{\hat{S}}\) respectively.
Recall that we have the suspension map \(\Sigma \co
{E}_{l-1}(\spc[k-1]{F}) \to {E}_l(\spc{F})\).  Let \(\Sigma_S = \pi_S
\Sigma \iota_S\) and \(\Sigma_{\hat{S}} = \pi_{\hat{S}} \Sigma
\iota_{\hat{S}}\).  Using these maps, we can consider the colimits of
the families \((s {E}_*(\spc{F}))\) and \(({E}_*(\spc{F})/ s
{E}_*(\spc{F}))\).

\begin{proposition}
The maps \(\iota_S\) etc.\ induce maps on the colimits.
\end{proposition}

\begin{proof}
To do this we need to show that they satisfy identities such as
\(\iota_S (\Sigma_S)^l = \Sigma^l \iota_S\) for some \(l\).  We shall
see that this works for the \(E\)\hyp{}unstable loop height of \({F}^*(-)\),
\(h\).   Since \(\iota_S \pi_S = S\),
\begin{align*}
\iota_S (\Sigma_S)^h &= \iota_S (\pi_S \Sigma \iota_S)^h \\
&= \iota_S \pi_S \Sigma \iota_S \pi_S \Sigma \iota_S \dotsb \pi_S
\Sigma \iota_S \\
&= S \Sigma S \Sigma \dotsb S \Sigma \iota_s \\
&= (S \Sigma)^h \iota_S.
\end{align*}
Now \(S\) is \(\circ\)\hyp{}multiplication by \(s\) and \(\Sigma\) is
(up to sign) \(\circ\)\hyp{}multiplication by \(e\).  As \(s \circ e = e
\circ s\) these two operations commute.  Furthermore, as \(s \circ
e^{\circ h} = e^{\circ h}\) these operations satisfy \(S \Sigma^h =
\Sigma^h\).  Putting this together yields the desired identity:
\[
  \iota_S (\Sigma_S)^h = S^h \Sigma^h \iota_S = \Sigma^h \iota_S.
\]
The other cases are similar; some use the fact that \(\pi_S \iota_S =
1\).
\end{proof}

\begin{corollary}
There is a split short exact sequence of algebras:

\hspace*{\fill}\(\displaystyle
  0 \to \colim_k s {E}_*(\spc{F}) \to \colim_k {E}_*(\spc{F}) \to \colim_k
  {E}_*(\spc{F}) / s {E}_*(\spc{F}) \to 0.\) \noproof
\end{corollary}

The colimits on the left and right are easily identified.

\begin{proposition}
The colimit on the right is null whereas the colimit of the left is
isomorphic to any of its components; that is, the natural map:
\[
  s {E}_i(\spc[j]{F}) \to \colim_k s {E}_{i+k}(\spc[j+k]{F})
\]
is an isomorphism.
\end{proposition}

\begin{proof}
This follows from the fact that \(\Sigma^h\) is an isomorphism on \(s
{E}_*(\spc[*]{F})\) and null on \({E}_*(\spc[*]{F})/ s {E}_*(\spc[*]{F})\),
where \(h\) is the \(E\)\hyp{}unstable loop height of \({F}^*(-)\).
\end{proof}

The ``null'' part implies that there is an isomorphism of algebras:
\[
  \colim s {E}_*(\spc[*]{F}) \cong \colim {E}_*(\spc[*]{F});
\]
whilst the other part implies that we can map back from the first
colimit to any of its components.

\begin{defn}
\label{def:desusp}
For \(k,l \in \Z\) let \(\delta \co \rc{E}_l(F) \to {E}_{l+k}(\spc{F})\)
be the map:
\[
  \rc{E}_l(F) \cong \colim_k {E}_{l+k}(\spc{F}) \cong
  \colim_k s {E}_{l+k}(\spc{F}) \cong s {E}_{l+k}(\spc{F}) \to
  {E}_{l+k}(\spc{F}),
\]
where the isomorphisms are as above.  We refer to \(\delta\) as the
\emph{destabilisation} map.
\end{defn}

\begin{proposition}
\label{prop:desusp}
The destabilisation map \(\delta\) is right-inverse to the
stabilisation map \({\sigma_k}_* \co {E}_{l+k}(\spc{F}) \to
\rc{E}_l(F)\).  The image of \(\delta\) is the image of the iterated
suspension map \(\Sigma^h \co {E}_{l+k-h}(\spc[k-h]{F}) \to
{E}_{l+k}(\spc{F})\) where \(h\) is the \(E\)\hyp{}unstable loop
height of \({F}^*(-)\).  In the particular case \(k = l = 0\),
\(\delta\) is a homomorphism of algebras. \noproof
\end{proposition}

The whole of the above can also be done in the additive realm and the
two correspond under the quotient map.

\section{Operations and Maps}
\label{sec:maps}

The results of the previous section readily dualise to operations due
to our assumption that operations from \({F}^*(-)\) to \({E}^*(-)\) are
dual to co-operations.  In this section we interpret our results
in the languages of operations and maps.  It will be obvious from this
formulation that the dual of the destabilisation map respects
composition of operations and maps.  We now state our main theorem.

\begin{theorem}
\label{th:opsmaps}
Let \({E}^*(-)\) and \({F}^*(-)\) be two graded multiplicative cohomology
theories that are commutative and complex orientable.  Suppose in
addition that the following conditions hold.
\begin{enumerate}
\item The coefficient ring, \(E^*\), of \({E}^*(-)\) has characteristic
\(p\).

\item The formal group law of \({E}^*(-)\) has finite height, say \(n\).

\item The coefficient of the first term in the \(p\)\hyp{}series for
\({E}^*(-)\) is invertible.

\item The various \(E^*\)\hyp{}modules of operations from \({F}^*(-)\)
to \({E}^*(-)\) are the \(E^*\)\hyp{}duals to the corresponding
\(E^*\)\hyp{}modules of co-operations.
\end{enumerate}
Let \(h\) be the \(E\)\hyp{}unstable loop height of \({F}^*(-)\).

Then there is a delooping map:
\begin{align*}
  \Delta^\infty &\co \m{U}^{k+l}_k(F,E) \to \m{S}^{l}(F,E) \\
\text{equivalently: } \Delta^\infty &\co {E}^{k+l}(\spc{F}) \to
\rc{E}^{l}(F) \\
\text{and: } \Delta^\infty &\co \lb \spc{F}, \spc[k+l]{E} \rb \to \{F, E\}^{l}
\end{align*}
left-inverse to the natural restriction map; thus in the last
formulation \(\Delta^\infty \Omega^\infty\) is the identity on \(\{F,
  E\}^{l}\).

Let \(r_k \in \m{U}^{k+l}_k(F,E)\) and let \(\rho_k \in
{E}^{k+l}(\spc{F})\) be the corresponding class.  The components of
\(\Delta^\infty r_k\) and \(\Delta^\infty \rho_k\) are:
\begin{align*}
(\Delta^\infty r_k)_m &= (-1)^{l m} (v_n^E)^{-i}(\Omega^j r_k)
(v_n^F)^i, \\
(\Delta^\infty \rho_k)_m &= (v_n^E)^{-i} (\Omega^j \rho_k)
(v_n^F)^i,
\end{align*}
where \(i,j \ge 0\) are chosen such that \(j \ge h\) and \(m - k =
2(p^n-1)i - j\).  In particular:
\begin{align*}
(\Delta^\infty r_k)_k &= (-1)^{l k} (v_n^E)^{-1} (\Omega^{2(p^n-1)}
r_k) v_n^F, \\
(\Delta^\infty \rho_k)_k &= (v_n^E)^{-1} (\Omega^{2(p^n-1)}
\rho_k) v_n^F; \\
\intertext{and for \(m \le k - h\):}
(\Delta^\infty r_k)_m &= (-1)^{l m} \Omega^{k-m} r_k, \\
(\Delta^\infty \rho_k)_m &= \Omega^{k-m} \rho_k.
\end{align*}

Moreover, an operation \(r_k \co {F}^k(-) \to {E}^{k+l}(-)\) is a
component of a stable operation if and only if it is the \(h\)\hyp{}fold
loop of an operation.  Similarly, a map \(\rho_k \co \spc{F} \to
\spc[k+l]{E}\) is an infinite loop map if and only if it is an
\(h\)\hyp{}fold loop map.
\end{theorem}

\begin{proof}
The delooping map is defined by dualising the destabilisation map,
\(\delta\), and using the correspondence between operations, maps,
and functionals to translate it across to the other realms.  As the
stabilisation map on co-operations is dual to the restriction map on
operations, the map \(\Delta^\infty\) is left-inverse to the natural
restriction map.

To determine the components of \(\Delta^\infty r_k\) and
\(\Delta^\infty \rho_k\) we first examine the components of an
arbitrary stable operation or map.  As \(h\) is the \(E\)\hyp{}unstable
loop height of \({F}^*(-)\), we have the identity:
\[
  v_n^E e^{\circ h} = e^{\circ (2(p^n-1) + h)} \circ \lb v_n^F \rb.
\]
Under stabilisation the element \(e\) maps to the identity
co-operation so the above stabilises to:
\[
  v_n^E = \lb v_n^F \rb.
\]
By assumption \(v_n^E\) is invertible.  Hence if \(c\) is a stable
co-operation \(c = (v_n^E)^{-1} c \lb v_n^F \rb\).  Dualising, if
\(r\) is a stable operation then \((v_n^E)^{-1} r \lb
v_n^F \rb = r\).  Let \((r_k)\) be the sequence of unstable operations
determined by restricting \(r\) to each degree.  The restriction maps
are bimodule maps and so we obtain the
identity:
\[
  r_k = (v_n^E)^{-1} r_{k - 2(p^n-1)} v_n^F.
\]
Now \(r_{k- 2(p^n-1)} = \Omega^{2(p^n-1)} r_k\) and hence:
\[
  r_k = (v_n^E)^{-1} (\Omega^{2(p^n-1)} r_k) v_n^F.
\]
Thus once we know one component of \(r\), say \(r_k\), we can
reconstruct the rest using the following procedure:
\begin{enumerate}
\item For \(m < k\) simply take the \((k-m)\)\hyp{}fold loop of \(r_k\).
\item For \(m > k\) take the \(j\)\hyp{}fold loop of \(r_k\) where
\(j\) is such that \(m-k+j = 2(p^n-1)i\) for some \(i > 0\).  Then the
periodicity ensures that:
\[
  r_m = (v_n^E)^{-i} r_{m - 2(p^n-1)i} (v_n^F)^i = (v_n^E)^{-i}
  r_{k-j} (v_n^F)^i = (v_n^E)^{-i} (\Omega^j r_k) (v_n^F)^i.
\]
\end{enumerate}

Thus the description of components of \(\Delta^\infty r_k\) and
\(\Delta^\infty \rho_k\) will follow from the final statement in the
theorem: that an \(h\)\hyp{}fold loop map is an infinite loop map.  This
is a direct consequence of the fact that the image of the
destabilisation map is the same as the image of the \(h\)th iterate of
the suspension map.  Hence the image of the delooping map is the image
of the \(h\)th iterate of the looping map.
\end{proof}

As \(\kn^*(-)\) has self unstable loop height of \(1\) we get the
following corollary.

\begin{corollary}
A map \(\alpha \co \spc{\kn} \to \spc[l]{\kn}\) is an infinite loop map
if and only if it is a loop map. \noproof
\end{corollary}

One further fact to record about the delooping map is that it respects
composition.

\begin{proposition}
Let \({E}^*(-), {F}^*(-), {G}^*(-)\) be graded multiplicative cohomology
theories such that the delooping maps \(\Delta^\infty_{F E}\),
\(\Delta^\infty_{G F}\), and \(\Delta^\infty_{G E}\) are all defined.
Let \(\rho_j \co \spc[j]{G} \to \spc{F}\) and \(\sigma_k \co \spc{F} \to
\spc[l]{E}\) be maps.  Then:
\[
  \Delta^\infty_{G E} (\sigma_k \rho_j) = \Delta^\infty_{F E}
  (\sigma_k) \Delta^\infty_{G F}(\rho_j).
\]
\end{proposition}

\begin{proof}
Firstly we note that both sides are well-defined.  Due to our
assumptions it is sufficient to show that this equation holds
component by component.  Moreover, due to the periodicity and the fact
that looping respects composition it is sufficient to show that it
holds for one component.  Thus we expand:
\begin{align*}
\big( \Delta^\infty_{F E}(\sigma_k) \Delta^{\infty}_{G F}(\rho_j)
\big)_j
&= \big( \Delta^\infty_{F E}(\sigma_k) \big)_k \big(\Delta^\infty_{G
  F} (\rho_j) \big)_j \\
&= (v_n^E)^{-1}(\Omega^{2(p^n-1)} \sigma_k) v_n^F (v_n^F)^{-1}
(\Omega^{2(p^n-1)} \rho_j) v_n^G \\
&= (v_n^E)^{-1}(\Omega^{2(p^n-1)} (\sigma_k \rho_j)) v_n^G;
\end{align*}
as required.
\end{proof}

\section{The Bousfield-Kuhn Functor}
\label{sec:bk}

In this section we relate our splitting to one that is a direct
consequence of the existence of the Bousfield-Kuhn functor.  In
\cite{nk2}, Kuhn showed that the \(\kn\)\hyp{}localisation of
\(p\)\hyp{}local spectra factors through the functor
\(\Omega^\infty\); this extended work of Bousfield in \cite{ab2} for
the case \(n = 1\).  For each \(n \ge 1\), Kuhn constructed a functor
\(\Phi_n\) from \(p\)\hyp{}local based spaces to \(p\)\hyp{}local
spectra such that \(\Phi_n \Omega^\infty\) is the
\(\kn\)\hyp{}localisation functor, \(L_{\kn}\).

As we now recall, the functorial properties of \(\Phi_n\) define a
map
\[
  \Theta_n \co \rkn^{k+l}(\spc{F}) \to \rkn^l(F)
\]
for any \(p\)\hyp{}local spectrum \(F\).  To see this, let \(\spc{F}\)
be the zeroth space of \(\Sigma^k F\) and recall that
\(\spc[k+l]{\kn}\) is the zeroth space of the \(p\)\hyp{}local
spectrum \(\Sigma^{k+l} \kn\).  As \(\Phi_n\) is a functor from
\(p\)\hyp{}local based spaces to \(\kn\)\hyp{}local spectra it defines
a map on morphism sets:
\[
  \lb \spc{F}, \spc[k+l]{\kn} \rb_+ \to \{L_{\kn} \Sigma^k F,
    L_{\kn} \Sigma^{k+l} \kn \}^0.
\]
We can simplify the target of this map.  The spectrum \(\Sigma^{k+l}
\kn\) is already \(\kn\)\hyp{}local allowing us to drop the second
\(L_{\kn}\).  This, together with sorting out the suspensions, means
that the target is naturally isomorphic to \(\{L_{\kn} F, \kn\}^l\)
which is \(\rkn^l(L_{\kn} F)\).  On the other hand, the source is
\(\rkn^{k+1}(\spc{F})\).  Hence we have a map:
\[
  \rkn^{k+l}(\spc{F}) \to \rkn^l(L_{\kn} F) \cong \rkn^l(F).
\]
In addition, we can remove the condition that \(F\) be \(p\)\hyp{}local by
noting that the Morava K\hyp{}theory is \(p\)\hyp{}local and hence only
``sees'' the \(p\)\hyp{}localisation of \(F\).  We can also extend
this to the unreduced cohomology theory using the canonical projection
of unreduced onto reduced cohomology.  Thus for any spectrum \(F\) the
Bousfield\hyp{}Kuhn functor defines a map
\[
  \Theta_n \co  \kn^{k+l}(\spc{F}) \to  \rkn^l(F).
\]

\begin{theorem}
Let \({F}^*(-)\) be a graded multiplicative cohomology theory that is
commutative and complex orientable.  Then the delooping map
\[
\Delta^\infty \co \kn^{k+l}(\spc{F}) \to \rkn^l(F)
\]
is defined and agrees with \(\Theta_n\).
\end{theorem}

\begin{proof}
The pair \(\kn^*(-)\) and \({F}^*(-)\) satisfy all the conditions
for the construction of \(\Delta^\infty\) and so it is at least
defined.  The map \(\Theta_n\) factors through the projection to
reduced cohomology by construction.  The same is true for
\(\Delta^\infty\) as can be seen from the formula in
theorem~\ref{th:opsmaps}: recall that our definition of the loop of an
unbased map involved first projecting it to a based map and then
taking the usual loop of the result.  Therefore to show that
\(\Delta^\infty\) and \(\Theta_n\) agree it is sufficient to show that
they agree on reduced cohomology; equivalently that they agree on
based maps.  In this situation the loop of a map is as expected with
no initial projection to based maps. 

The first step in showing that \(\Delta^\infty\) and \(\Theta_n\) are
the same map is to observe that, as both are left-inverse to
\(\Omega^\infty\), if a class \(\rho_k\) is a component of a stable
class then \(\Delta^\infty (\rho_k) = \Theta_n (\rho_k)\).  By
theorem~\ref{th:opsmaps} any unstable class becomes the component of a
stable class after a finite number of loopings.  Therefore it is
enough to show that \(\Delta^\infty\) and \(\Theta_n\) both commute
with loops.  In both cases this is immediate from the constructions of
the maps.  For completeness we review the definition of the
Bousfield-Kuhn functor from \cite{nk2} and explain how the desired
property follows.

There are three steps in defining \(\Phi_n\).
\begin{enumerate}
\item Let \(Z\) be a finite CW-complex with a self-map \(v \co \Sigma^d Z
\to Z\), \(d > 0\).   Composition with \(v\) defines a map
\[
  v^* \co \map(Z,X) \to \map(\Sigma^d Z, X) = \Omega^d \map(Z,X)
\]
for any based space \(X\).  One can therefore define a spectrum with
\(m d\)th space \(\map(Z,X)\) and structure maps \(v^*\).  This
construction is functorial in \(X\) and so defines a functor
\(\Phi_Z'\) from based spaces to spectra.

\item The second step is to compose the functor \(\Phi_Z'\) with
\(\kn\)\hyp{}localisation to produce a functor \(\Phi_Z\) from spaces to
\(\kn\)\hyp{}local spectra.

\item The final step is to define a functor \(\Phi_n\) from spaces to
\(\kn\)\hyp{}local spectra by taking the direct limit of a sequence of
functors, \((\Phi_{Z_k})\), for a suitable choice of sequence of
spaces \((Z_k)\).
\end{enumerate}

Both localisation and taking the direct limit of a sequence of spectra
commute with the suspension and loop operators acting on the category
of spectra.  Therefore to show that \(\Phi_n\), and thus \(\Theta_n\),
commutes with looping it is sufficient to show that this is true for
\(\Phi_Z'\).  This follows from the fact that \(\map(Z, \Omega X) =
\Omega \map(Z, X)\).  Thus the spectrum for \(\Phi_Z'(\Omega X)\) is
the spectrum \(\Omega \Phi_Z'(X)\) and similarly \(\Phi_Z'(\Omega
\alpha) = \Omega \Phi_Z'(\alpha)\) for a based map \(\alpha \co X \to
Y\).
\end{proof}

\end{document}